\documentclass[letterpaper, 10 pt, conference]{ieeeconf}
\IEEEoverridecommandlockouts                              
\usepackage{longtable}
\usepackage{booktabs}

\newif\ifdeuxcols
\makeatletter
\if@twocolumn
  \deuxcolstrue
\else
  \deuxcolsfalse
\fi
\makeatother

\makeatletter
\newcommand{\purpose}[1]{\def\@purpose{#1}}
\purpose{}
\def\@oddfootlefttext{
  \ifx\@purpose\@empty
    Preprint submitted to \ifx\@journal\@empty Elsevier\else\@journal\fi
  \else\@purpose\fi}

\def\ps@pprintTitle{%
     \let\@oddhead\@empty
     \let\@evenhead\@empty
     \def\@oddfoot{\footnotesize\itshape
       \@oddfootlefttext\hfill\today}%
     \let\@evenfoot\@oddfoot}
\makeatother

\usepackage{amsmath} 
\usepackage{epsfig}
\usepackage{graphicx}
\usepackage[hang]{caption2}
\usepackage[normalsize,it,hang]{subfigure}
\usepackage[latin1]{inputenc}
\usepackage[greek,english]{babel}
\newtheorem{remark}{Remark}[section]

\title{\LARGE \bf
A mathematical explanation via ``intelligent'' PID controllers \\ of
the strange ubiquity of PIDs }

\author{Brigitte d'{\sc Andr\'{e}a-Novel},  Michel {\sc Fliess},  C\'edric {\sc Join},
Hugues {\sc Mounier},  Bruno {\sc Steux}
\thanks{Brigitte d'{\sc Andr\'{e}a-Novel} is with Centre de Robotique, Mines-ParisTech,
75272 Paris Cedex 06, France.
\newline        {\tt\small Brigitte.Dandrea-Novel@mines-paristech.fr}}
\thanks{Michel {\sc Fliess} is with INRIA-ALIEN \& LIX (CNRS, UMR 7161),
\'Ecole polytechnique, 91128 Palaiseau, France.
\newline        {\tt\small Michel.Fliess@polytechnique.edu}}
\thanks{C\'{e}dric {\sc Join} is with INRIA-ALIEN \& CRAN (CNRS, UMR 7039), Nancy-Universit\'e,
BP 239, 54506 Vand\oe uvre-l\`es-Nancy, France.
\newline        {\tt\small cedric.join@cran.uhp-nancy.fr}}
\thanks{Hugues {\sc Mounier} is with L2S (CNRS, UMR 8506),
 Sup\'{e}lec \& Universit\'{e} Paris-Sud, 3 rue Joliot-Curie, 91192 Gif-sur-Yvette, France.
\newline        {\tt\small Hugues.Mounier@lss.supelec.fr}}
\thanks{Bruno {\sc Steux} is with Centre de Robotique, Mines-ParisTech,
75272 Paris Cedex 06, France.
\newline {\tt\small Bruno.Steux@mines-paristech.fr}}
        }

\begin{document}

\maketitle
\thispagestyle{empty}
\pagestyle{empty}






\begin{abstract}
The ubiquity of PID controllers in the industry has remained
mysterious until now. We provide here a mathematical explanation of
this strange phenomenon by comparing their sampling with the
the one of ``intelligent'' PID controllers, which were recently
introduced. Some computer simulations  nevertheless confirm the
superiority of the new intelligent feedback design.
\end{abstract}


\section{Introduction}
\label{intro} PI and PID controllers (see, {\it e.g.},
\cite{pid,od}) are still by far the most popular feedback design in
industry. To the best of our knowledge, there is no clear-cut
explanation of their strange ubiquity for a wide range of systems.
Remember that, from a purely mathematical standpoint, they are only
fully justified until now for first and second order linear
differential equations with constant coefficients! We solve here
this long-standing and quite irritating open problem via the newly
introduced {\em intelligent} PIDs (\cite{esta,malo}), which have
already been utilized quite successfully in several concrete
situations (see, {\it e.g.},
\cite{cifa-mines,choi,brest1,brest2,psa,edf,vil}).

The proof relies on a crude time-sampling of both types of
regulators. It shows that the gains in a classic PI or PID take into
account, if they are properly tuned, the estimated ``structural''
part of the intelligent controllers. Thus the efficiency of these
intelligent controllers with respect to arbitrary nonlinear plants
(\cite{esta,malo}) is enough for fulfilling our purpose. Let us
nevertheless emphasize that the classic tuning rules are quite
intricate whereas their counterparts for intelligent controllers are
obvious.

\begin{remark}
Only few references (see, {\it e.g.}, \cite{chang-jung,han}) in the
huge literature on PIDs exhibit some connections with our
viewpoint.
\end{remark}

Our paper is organized as follows. Section \ref{mfc} is devoted to a
brief review of {\em model-free} control and of the corresponding
intelligent PID controllers. Section \ref{connect}, which
establishes our results by comparing the sampling of classical and
intelligent controllers, gives a table for the connections between
classic and intelligent gains. The computer simulations in Section
\ref{Applications} confirms the superiority of the intelligent
controllers (see, also, \cite{esta,malo} for other examples). Some
concluding remarks are given in Section \ref{Conclusion}.

\section{Model-free control\protect\footnote{See
\cite{esta,malo} for more details.}}\label{mfc}

The input-output behavior of the system, which for simplicity's
sake is assumed to be
monovariable, is ``approximatively'' governed within its operating
range by an {\bf unknown} finite-dimensional ordinary differential
equation, which is not necessarily linear,
\begin{equation}\label{E}
E (y, \dot{y}, \dots, y^{(a)}, u, \dot{u}, \dots, u^{(b)}) = 0
\end{equation}
Replace Equation \eqref{E} by the following ``pheno\-menological''
model, which is only valid during a very short time interval,
\begin{equation}\label{F}
y^{(\nu)} = F + \alpha u
\end{equation}
The derivation order $\nu$, which is in general equal to $1$ or $2$,
and the constant parameter $\alpha$ are chosen by the practitioner.
It implies that $\nu$ is not necessarily equal to the derivation
order $a$ of $y$ in Equation (\ref{E}). The numerical value of $F$
at any time instant is deduced from those of $u$ and $y^{(\nu)}$,
thanks to our quite efficient numerical differentiators, which are
moreover real-time.\footnote{See \cite{easy}, \cite{mamadou} for
details, and also \cite{garcia}.} The desired behavior is obtained
by implementing, if, for instance, $\nu = 2$, the {\em intelligent
PID controller} ({\em i-PID})
\begin{equation}\label{universal}
u = \frac{1}{\alpha} \left(- F + \ddot{y}^\ast + K_P e + K_I \int e
+ K_D \dot{e} \right)
\end{equation}
where
\begin{itemize}
\item $y^\ast$ is the output reference trajectory, which
is determined e.g. via the rules of flatness-based control;
\item $e = y - y^\ast$ is the tracking error;
\item $K_P$, $K_I$, $K_D$ are the usual tuning gains.
\end{itemize}

Let us consider the following special cases:
\begin{itemize}
\item If again $\nu = 2$, we may use an {\em intelligent PD controller} ({\em i-PD})
\begin{equation}\label{iPD}
u = \frac{1}{\alpha} \left(-F + \ddot{y}^\ast + K_P e + K_D \dot{e}
\right)
\end{equation}
\item If $\nu = 1$, we can restrict ourselves to
\begin{itemize}
\item an {\em intelligent PI controller} ({\em i-PI})
\begin{equation}\label{iPI}
u = \frac{1}{\alpha} \left(-F + \dot{y}^\ast + K_P e + K_I \int e
\right)
\end{equation}
\item or even to an {\em intelligent P controller} ({\em i-P})
\begin{equation}\label{iP}
u = \frac{1}{\alpha} \left(-F + \dot{y}^\ast + K_P e \right)
\end{equation}
\end{itemize}
\end{itemize}

\begin{remark}
If $\nu = 2$ (resp. $1$), plugging Equations \eqref{universal} or
\eqref{iPD} (resp. \eqref{iPI} or \eqref{iP}) in Equation \eqref{F}
yields the control of a pure double (resp. simple) integrator. This
is why tuning the gains of our intelligent controllers is quite
straightforward.
\end{remark}

\begin{remark}
It should be emphasized, if $\nu = 2$ (resp. $1$), that Equation
\eqref{iPD} (resp. \eqref{iP}) is mathematically sufficient for
ensuring stability around the reference trajectory. The integral
term $ K_I \int e$ in Equation \eqref{universal} (resp. \eqref{iPI})
nevertheless adds some well known robustness properties.
\end{remark}

\section{Connections between classic and intelligent
controllers}\label{connect}
\subsection{PI and i-P}\label{1}
\subsubsection{A crude sampling of PIs}
%
Consider the classic continuous-time PI controller
\begin{equation}\label{cpi}
  u (t) = k_p e(t) + k_i \int e(\tau) d\tau
\end{equation}
A crude sampling of the integral $\int e(\tau) d\tau$ through a
Riemann sum $I(t)$ leads to
$$
\int e(\tau) d\tau \simeq  I(t) = I(t-h) + h e(t)
$$
where $h$ is the sampling interval. The corresponding discrete form
of Equation \eqref{cpi} reads:
$$
u(t) = k_p e(t) + k_i I(t) = k_p e(t) + k_i I(t-h) + k_i h e(t)
$$
Combining the above equation with $$u(t-h) = k_p e(t-h) + k_i
I(t-h)$$ yields
\begin{equation}
\label{eqPIRiemannDiscrSix} u(t) = u(t - h) + k_p \left( e(t) - e(t
- h) \right) + k_i h e(t)
\end{equation}

\begin{remark}
A trivial sampling of the ``velocity form'' of Equation \eqref{cpi}
\begin{equation*}\label{cpid}
  \dot{u} (t) = k_p \dot{e}(t) + k_i e(t)
\end{equation*}
yields
$$
\dfrac{u(t) - u(t - h)}{h} =  k_p  \left(\dfrac{e(t) - e(t -
h)}{h}\right) + k_i  e(t)
$$
which is equivalent to Equation \eqref{eqPIRiemannDiscrSix}.
\end{remark}

\subsubsection{Sampling i-Ps} Utilize, if $\nu = 1$ in
Equation \eqref{F}, the i-P \eqref{iP}, which may be rewritten as
$$
u (t) = \frac{{\dot y}^\ast (t) - F + K_P e(t)}{\alpha}
$$
Replace, according to the computer implementation in
\cite{esta,malo}, $F$ by ${\dot y}(t) - \alpha u (t-h)$ and
therefore by
$$\frac{y(t) - y(t-h)}{h} - \alpha u (t-h)$$
It yields
\begin{equation}
\label{eqDiscr_i-POne} u (t) = u (t - h) - \frac{e(t) -
e(t-h)}{h\alpha} + \dfrac{K_P}{\alpha}\, e(t)
\end{equation}
\subsubsection{Comparison}
{\bf FACT}.- Equations \eqref{eqPIRiemannDiscrSix} and
\eqref{eqDiscr_i-POne} become {\bf identical} if we set
\begin{align}
\label{eqPI_i-P_corresp} k_p &= - \dfrac{1}{\alpha h}, \quad k_i =
\dfrac{K_P}{\alpha h}
\end{align}

\begin{remark}
It should be emphasized that the above property, defined by
Equations \eqref{eqPI_i-P_corresp}, does not hold for
continuous-time PIs and i-Ps. This equivalence is strictly related
to time sampling, {\em i.e.}, to computer implementation, as
demonstrated by taking $h \rightarrow 0$ in Equations
\eqref{eqPI_i-P_corresp}.
\end{remark}

\subsection{PID and i-PD}
Extending the calculations of Section \ref{1} is quite obvious. The
velocity form of the PID
$$u(t) = k_p e(t) + k_i \int e(\tau) d\tau + k_d \dot{e}$$ reads
$\dot u(t) = k_p {\dot e}(t) + k_i e(t) + k_d {\ddot e}(t)$. It
yields the obvious sampling
\begin{equation}
\label{eqDiscr_PID}
  u(t) = u(t - h) + k_p h {\dot e} (t) + k_i h e(t) +
            k_d h {\ddot e}(t)
\end{equation}
If $\nu = 2$ on the other hand, Equation \eqref{iPD} yields $u (t) =
\dfrac{1}{\alpha} \left( {\ddot y}^\ast (t) - F + K_P e(t) + K_D
{\dot e}(t) \right)$. From the computer implementation $F =
\ddot{y}(t) - \alpha u(t - h)$, we derive
\begin{equation}\label{iPDd}
u (t) = u (t - h) - \dfrac{1}{\alpha} {\ddot e}(t)  +
\dfrac{K_P}{\alpha} e(t) + \dfrac{K_D}{\alpha} {\dot e}(t)
\end{equation}

\noindent{\bf FACT}.- Equations \eqref{eqDiscr_PID} and \eqref{iPDd}
become {\bf identical} if we set
\begin{equation}
\label{eqPID_i-PD_corresp} k_p = \dfrac{K_D}{\alpha h}, \quad k_i =
\dfrac{K_P}{\alpha h}, \quad k_d = - \dfrac{1}{\alpha h}
\end{equation}

\subsection{i-PI and i-PID}
\label{Seci-PIsi-PIDs} Equation \eqref{iPDd} becomes with the i-PID
\eqref{universal}
\begin{equation}\label{iPIDd}
u (t) = u (t - h) - \dfrac{1}{\alpha} {\ddot e}(t)  +
\dfrac{K_P}{\alpha} e(t) + \dfrac{K_I}{\alpha} \int e +
\dfrac{K_D}{\alpha} {\dot e}(t)
\end{equation}
Introduce the PII$^2$D controller
$$u (t) = k_p e (t) + k_i \int e(\tau)d\tau  + k_{ii}
\int\!\!\!\!\int e d\tau d\sigma + k_d \dot{e}(t)$$ where a double
integral appears.\footnote{Such double integrals do not seem to be
common in control engineering.} To its velocity form ${\dot u} (t) =
k_p \dot{e} (t) + k_i e + k_{ii} \int e d\tau + k_d \ddot{e}(t)$
corresponds the sampling
$$
u (t) = u(t-h) + k_p h \dot{e} (t) + k_i h e + k_{ii} h \int e d\tau
+ k_d h \ddot{e}(t) $$ which is identical to Equation \eqref{iPIDd}
if one sets
\begin{align}
    \label{eqPID_i-PD_correspTwo}
    k_p &= \dfrac{K_D}{\alpha h}, \quad k_i = \dfrac{K_P}{\alpha h}, \quad
    k_{ii} = \dfrac{K_I}{\alpha h}, \quad k_d = - \dfrac{1}{\alpha h}
  \end{align}
The connection between iPIs and PII$^2$s follows at once.
\subsection{Table of correspondence} \label{SecCorrespondance} The
previous calculations yield the following correspondence table
between the gains of our various controllers: \ifdeuxcols
\begin{footnotesize}
\fi
\begin{table*}[htb]\label{table}
\begin{center}
\begin{tabular}{p{5ex}p{3ex}p{7ex}p{7ex}p{7ex}p{7ex}}
\toprule
         &       & i-P            & i-PD           & i-PI & i-PID          \\
\midrule
PI       & $k_p$ & $- 1/\alpha h$   &                &      &                \\
         & $k_i$ & $K_P/\alpha h$ &                &      &                \\
\midrule
PID      & $k_p$ &                & $K_D/\alpha h$   &      &                \\
         & $k_i$ &                & $K_P/\alpha h$ &      &                \\
         & $k_d$ &                & $- 1/\alpha h$ &      &                \\
\midrule
PII$^2$  & $k_p$    &             &                & $- 1/\alpha h$   &      \\
         & $k_i$    &             &                & $K_P/\alpha h$ &      \\
         & $k_{ii}$ &             &                & $K_I/\alpha h$ &      \\
\midrule
PII$^2$D & $k_p$    &             &                &      & $K_D/\alpha h$   \\
         & $k_i$    &             &                &      & $K_P/\alpha h$ \\
         & $k_{ii}$ &             &                &      & $K_I/\alpha h$ \\
         & $k_d$    &             &                &      & $- 1/\alpha h$ \\
\bottomrule
\end{tabular}
\end{center}
\caption{\label{tblGainsCorresp}Correspondence between the gains of
sampled classic and intelligent controllers.}
\end{table*}
\ifdeuxcols
\end{footnotesize}
\fi

\begin{remark}
Due to the form of Equation (\ref{F}), it should be noticed that the
tuning gains of the classic regulators ought to be negative.
\end{remark}

\subsection{The explanation}
The previous calculations and Table 1 explain why sampled classic PI
and PID controllers take into account, if their gains are properly
tuned, the structural term $-F/\alpha$, which contains all the
structural information of the unknown nonlinear systems, in
Equations \eqref{universal}, \eqref{iPD}, \eqref{iPI}, \eqref{iP}.
The superiority of intelligent controllers, which was already noted
in \cite{esta,malo}, is however confirmed:
\begin{enumerate}
\item Tuning the gains of intelligent controllers is str\-aightforward
whereas it is complex and painful for classic PIDs in spite of all
the numerous existing rules in the literature (see, {\it e.g.},
\cite{pid,od}).
\item Contrarily to intelligent controllers, a correctly
tuned classic PI or PID controller is unable to take into account
heat effects, ageing processes, characteristic dispersions due to
mass production, \dots.
\item Fault tolerant control is much better handled by intelligent
controllers than by classic ones.
\end{enumerate}

\section{Classic versus intelligent controllers\protect\footnote{See \cite{esta,malo} for other
examples.}} \label{Applications} 
\begin{figure*}[htb]
\subfigure[\scriptsize
Input]{\includegraphics[scale=.405]{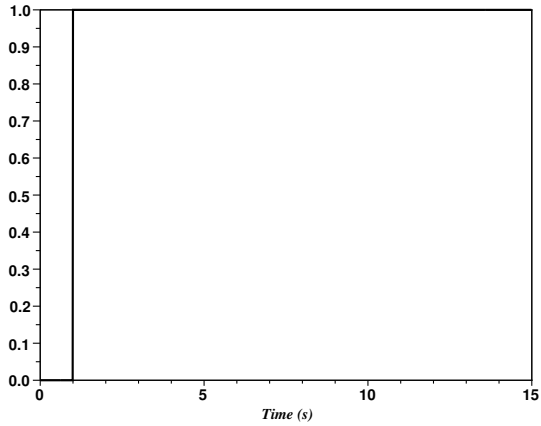}}
\subfigure[\scriptsize Output (--) and denoised output (-
-)]{\includegraphics[scale=.405]{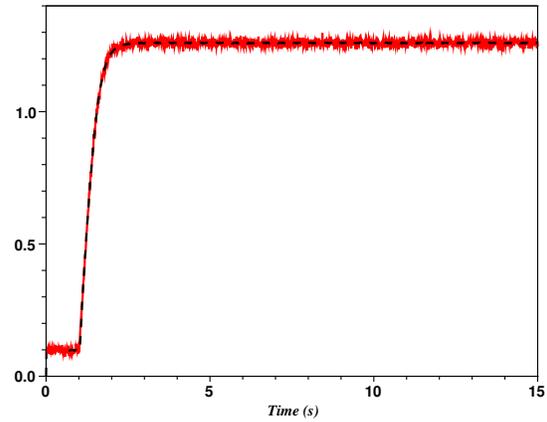}}
\caption{Open loop case}\label{x1BO}
\end{figure*}

\begin{figure*}[htb]
\subfigure[\scriptsize
Input]{\includegraphics[scale=.405]{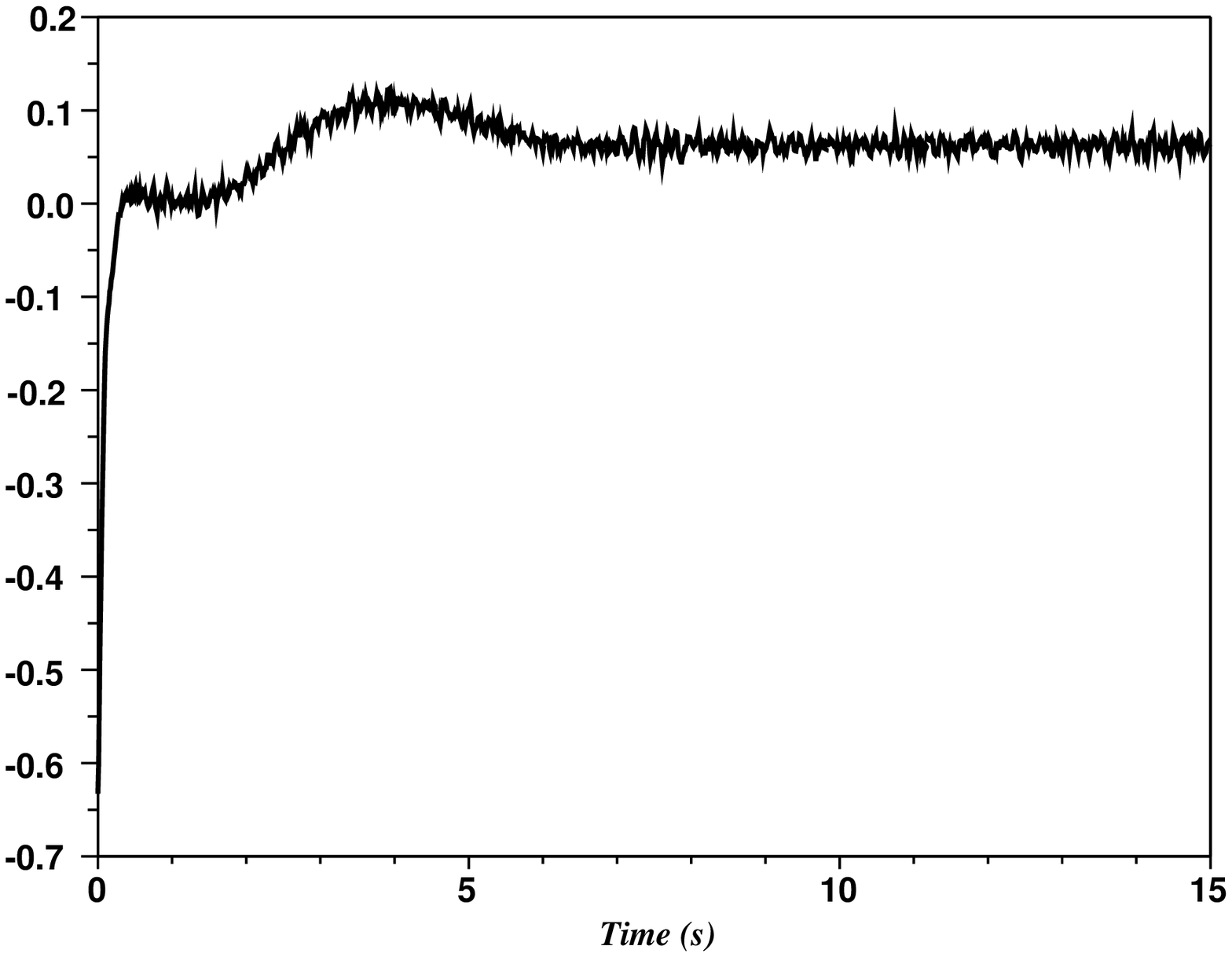}}
\subfigure[\scriptsize Output (--) and denoised output (-
-)]{\includegraphics[scale=.405]{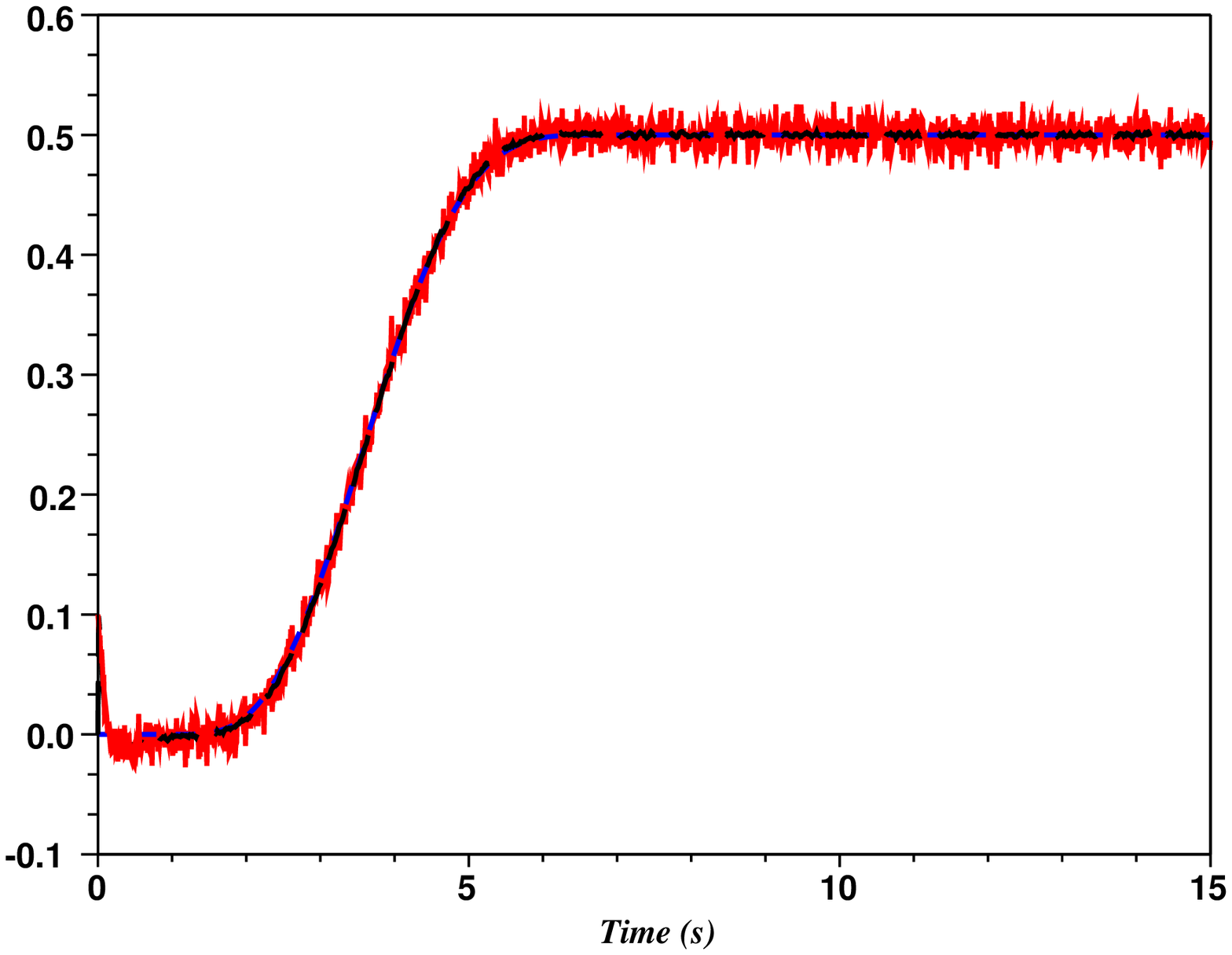}}
\caption{PI case}\label{x1PI}
\end{figure*}

\begin{figure*}[htb]
\subfigure[\scriptsize
Input]{\includegraphics[scale=.405]{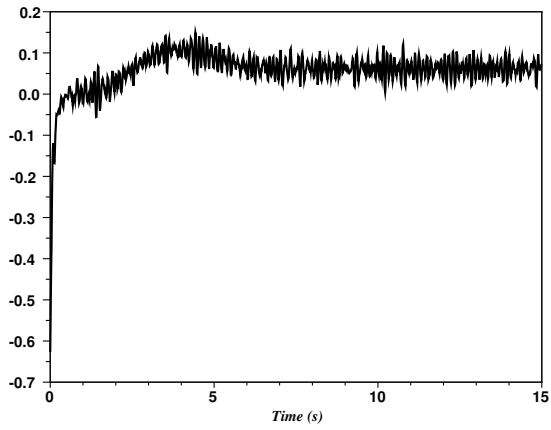}}
\subfigure[\scriptsize Output (--) and denoised output (-
-)]{\includegraphics[scale=.405]{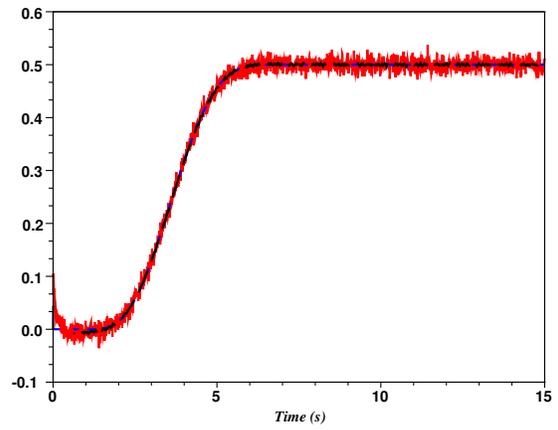}}
\caption{i-PI case}\label{x1iPI}
\end{figure*}

For the nonlinear system
\begin{equation}\label{nl}
\dot{y} +y^3=2u
\end{equation}
we deduce a classic PI controller thanks to a method due to
Bro\"{\i}da and Dindeleux \cite{dindeleux} which improves the
well-known Ziegler-Nichols rules (see, {\it e.g.}, \cite{pid,od}).
Note however that the open loop response of System \eqref{nl}, with
$y(0) = 0$, is somehow difficult to exploit as shown by Figure
\ref{x1BO}. It yields
\begin{itemize}
\item a delay system
$$\frac{ke^{-\tau s}}{1+Ts}$$
where $k=1.160$, $T=0.401$, $\tau=0.044$;
\item a PI where $k_p=6.350$, $k_i=15.817$.
\end{itemize}
Figures \ref{x1PI} and \ref{x1iPI}, which depict the simulation
results for the above PI and an i-PI, do not show any significant
difference. Remember however that the i-PI, where $\alpha=1$, $K_P
=6$, $K_I = 9$, does not necessitate any cumbersome identification
procedure.

\begin{figure*}[htb]
\subfigure[\scriptsize
Input]{\includegraphics[scale=.405]{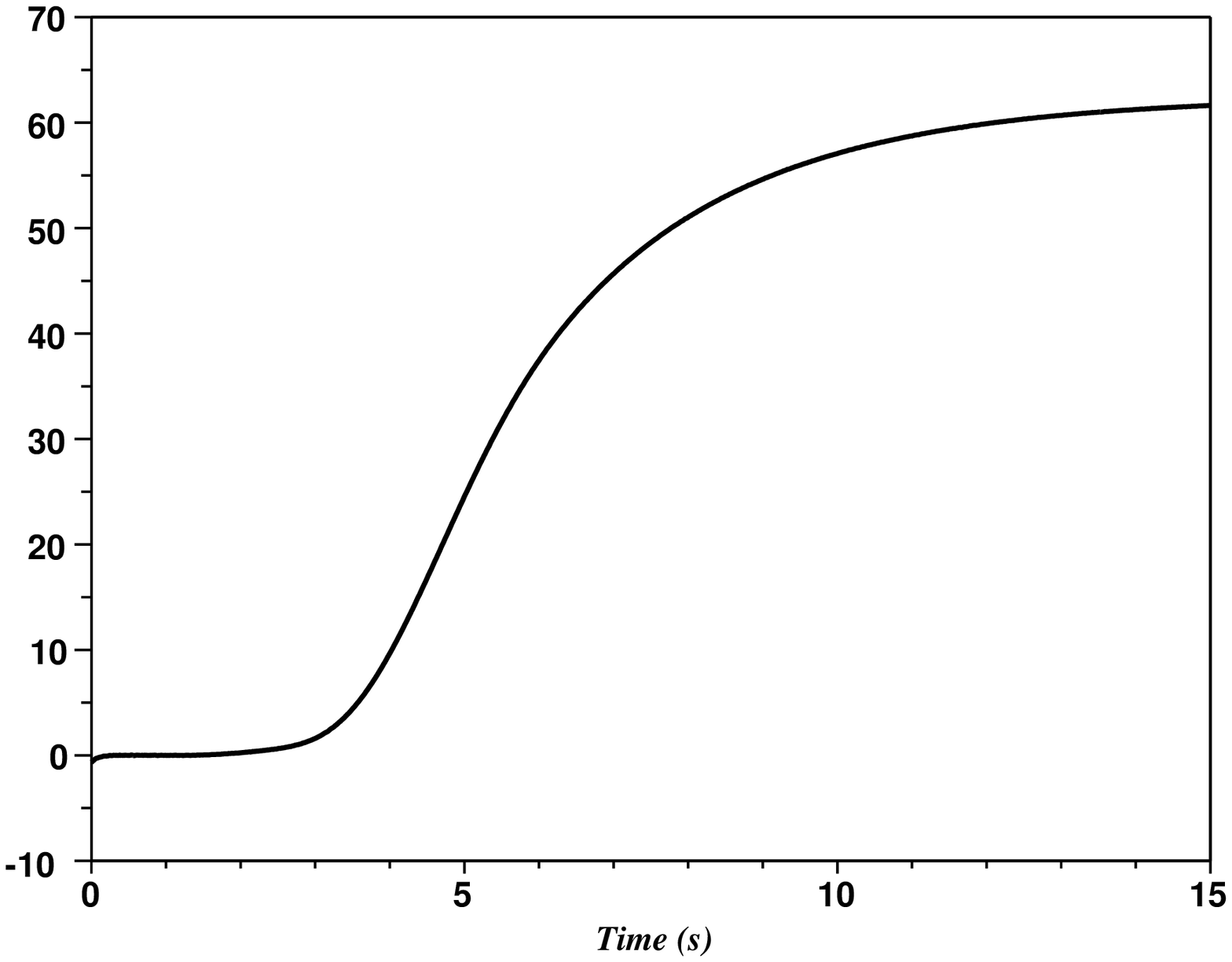}}
\subfigure[\scriptsize Output (--) and denoised output (-
-)]{\includegraphics[scale=.405]{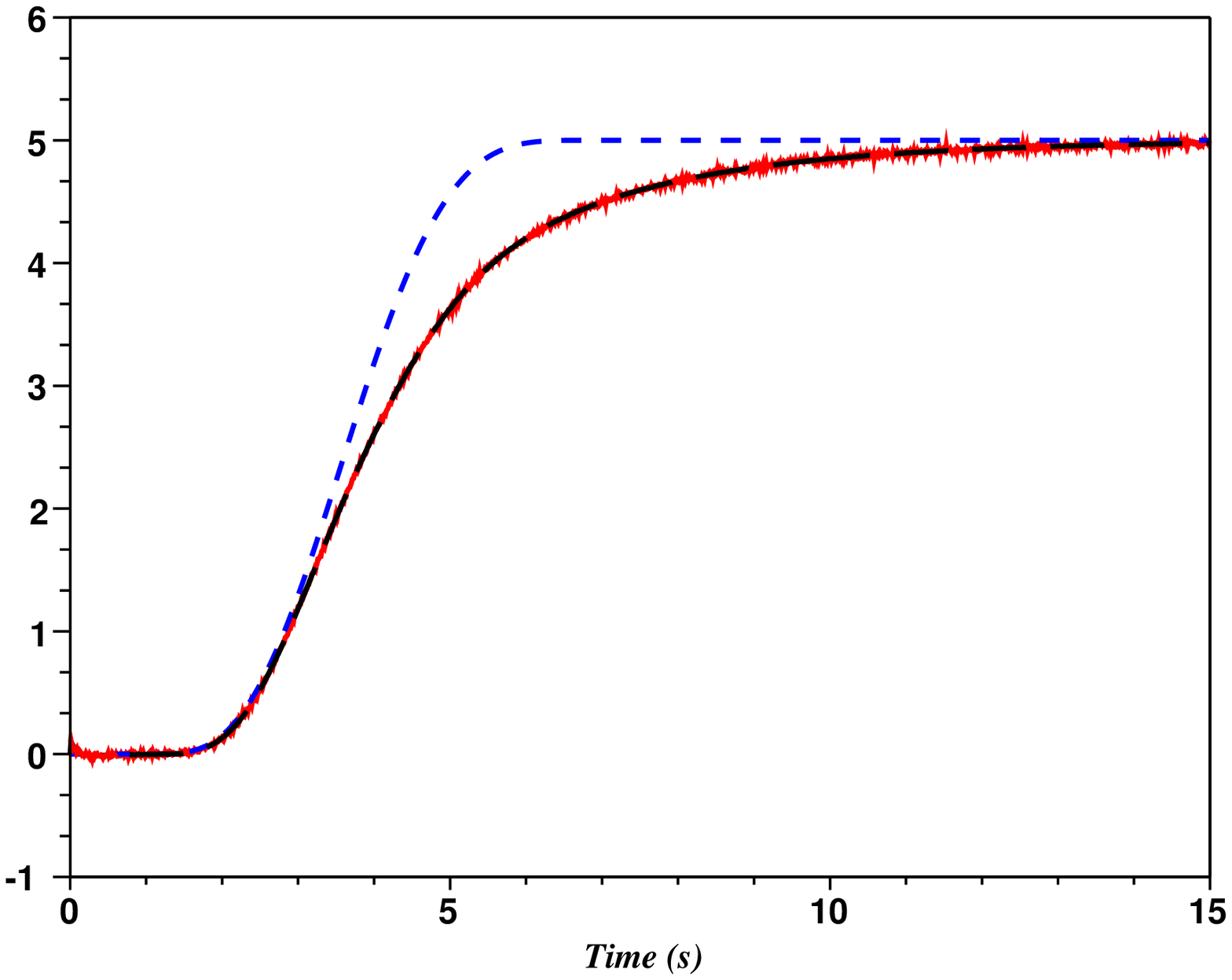}}
\caption{PI in case of large amplitude setpoint change}\label{x1PI5}
\end{figure*}

\begin{figure*}[htb]
\subfigure[\scriptsize
Input]{\includegraphics[scale=.405]{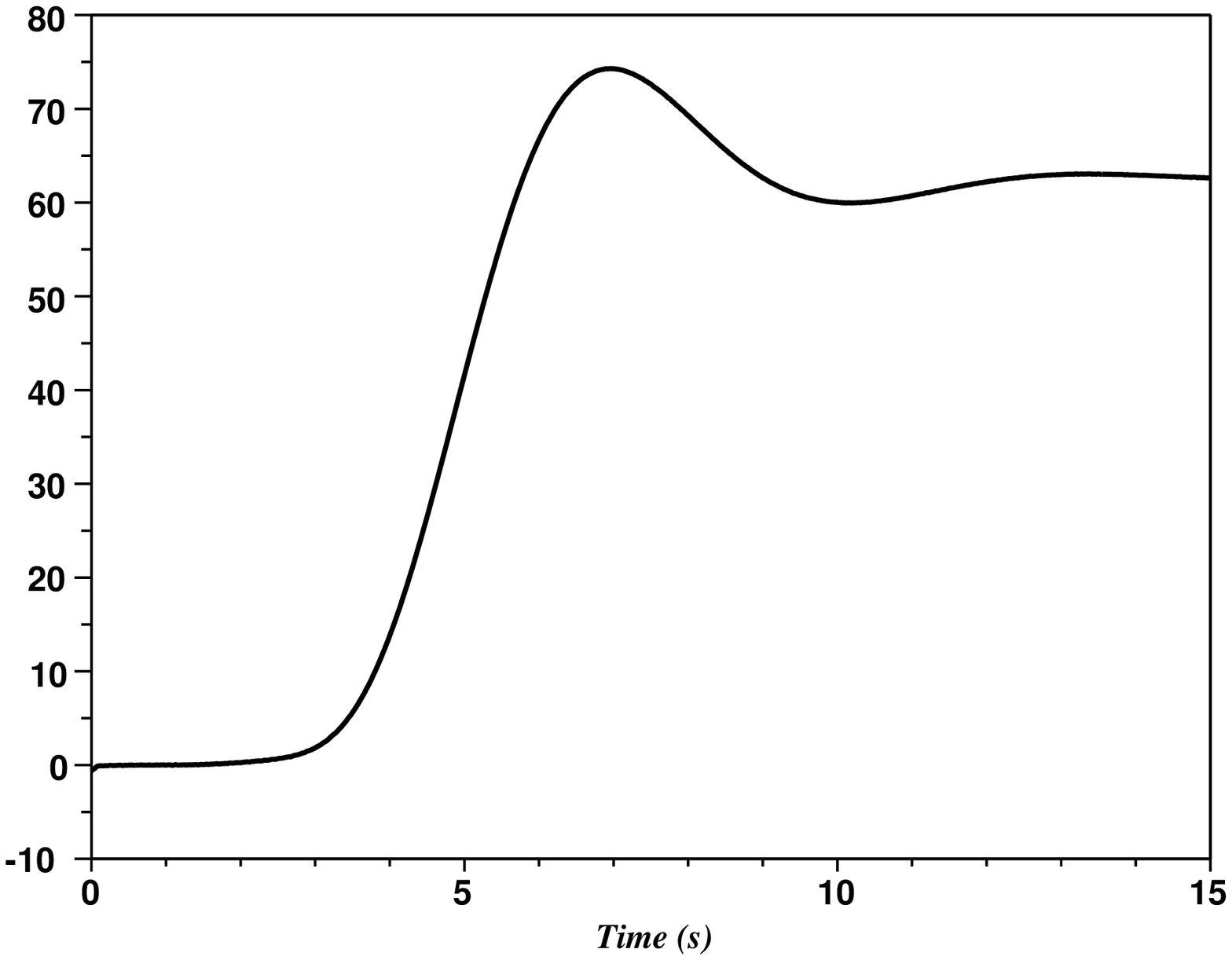}}
\subfigure[\scriptsize Output (--) and denoised output (-
-)]{\includegraphics[scale=.405]{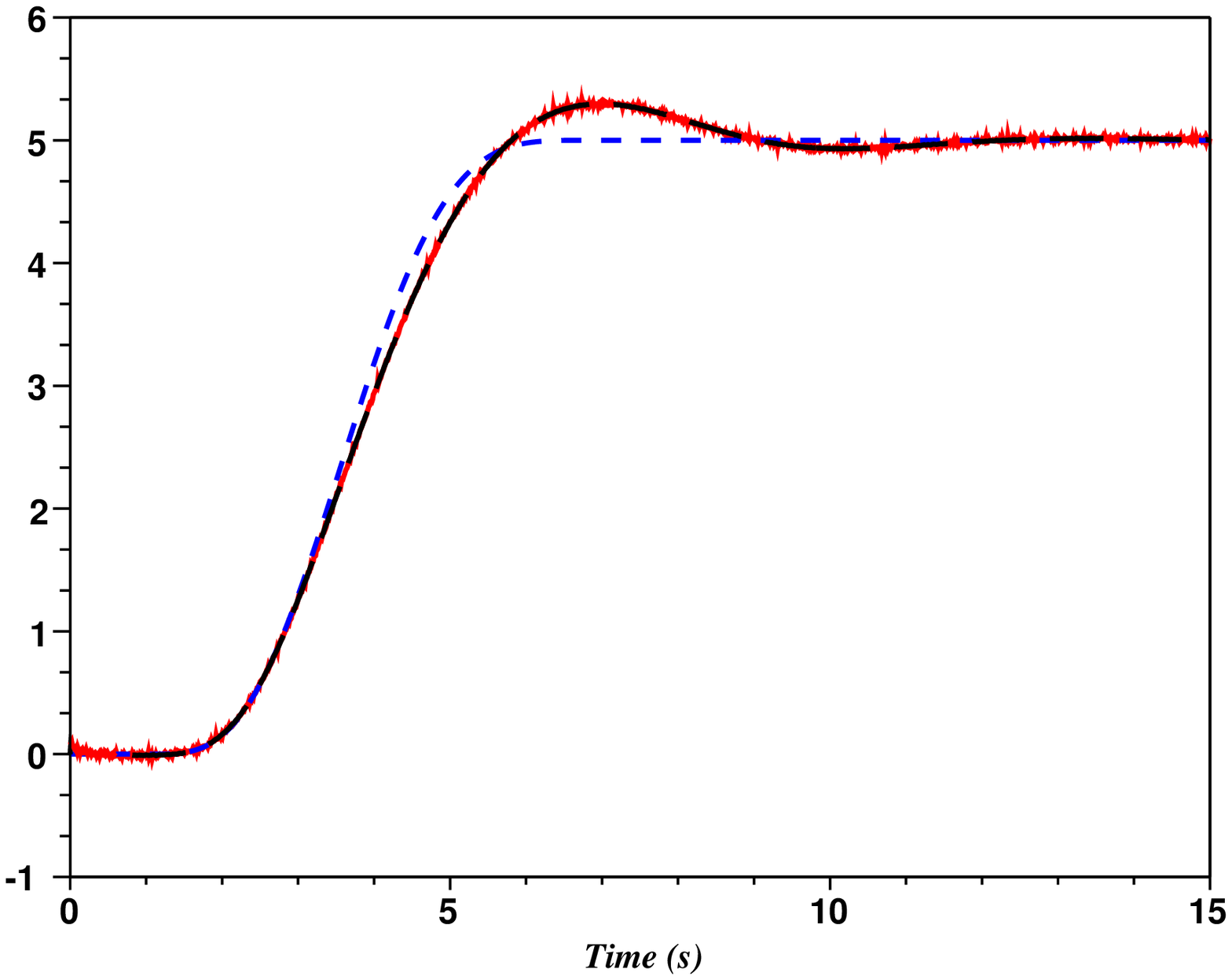}}
\caption{i-PI in case of large amplitude setpoint change}\label{x1iPI5}
\end{figure*}

Without any new calibration of the PI for another operating range
Figure \ref{x1PI5} shows a deterioration of the performances,
whereas the performances of the i-PI, which are depicted in Figure
\ref{x1iPI5}, remain good.

\begin{figure*}[htb]
\subfigure[\scriptsize
Input]{\includegraphics[scale=.405]{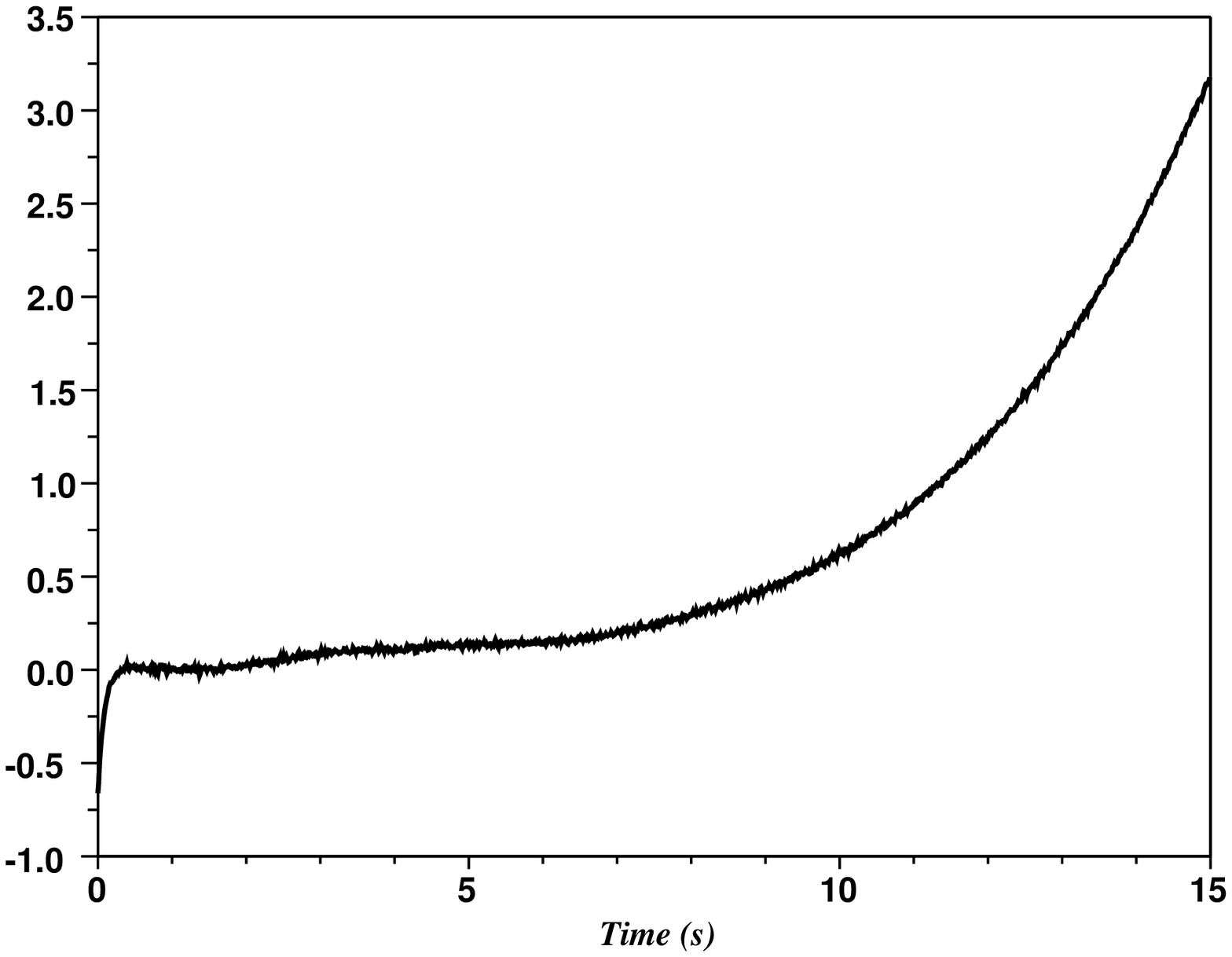}}
\subfigure[\scriptsize Output (--) and denoised output (-
-)]{\includegraphics[scale=.405]{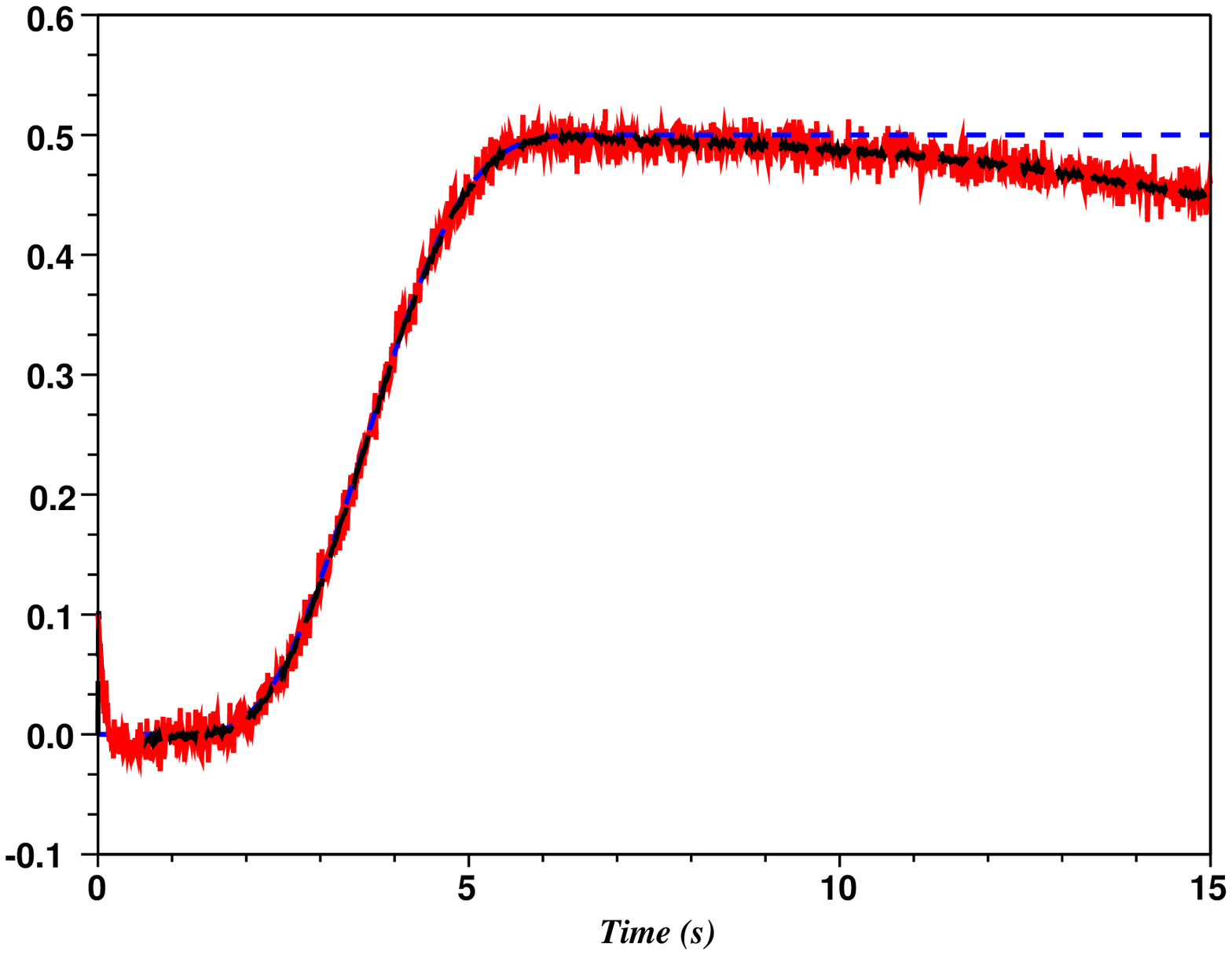}}
\caption{PI in case of power loss}\label{x1PIPP}
\end{figure*}

\begin{figure*}[htb]
\subfigure[\scriptsize
Input]{\includegraphics[scale=.405]{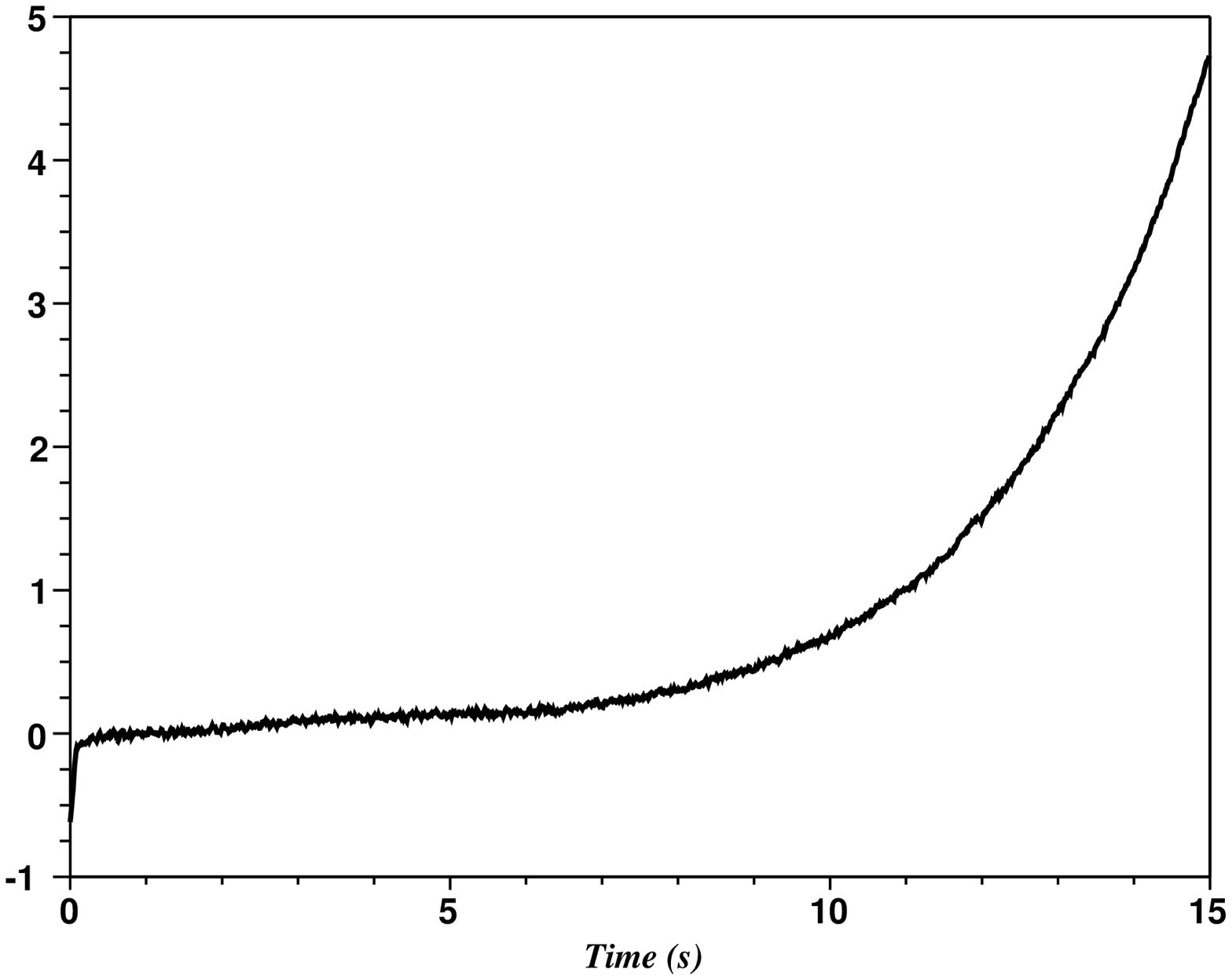}}
\subfigure[\scriptsize Output (--) and denoised output (-
-)]{\includegraphics[scale=.405]{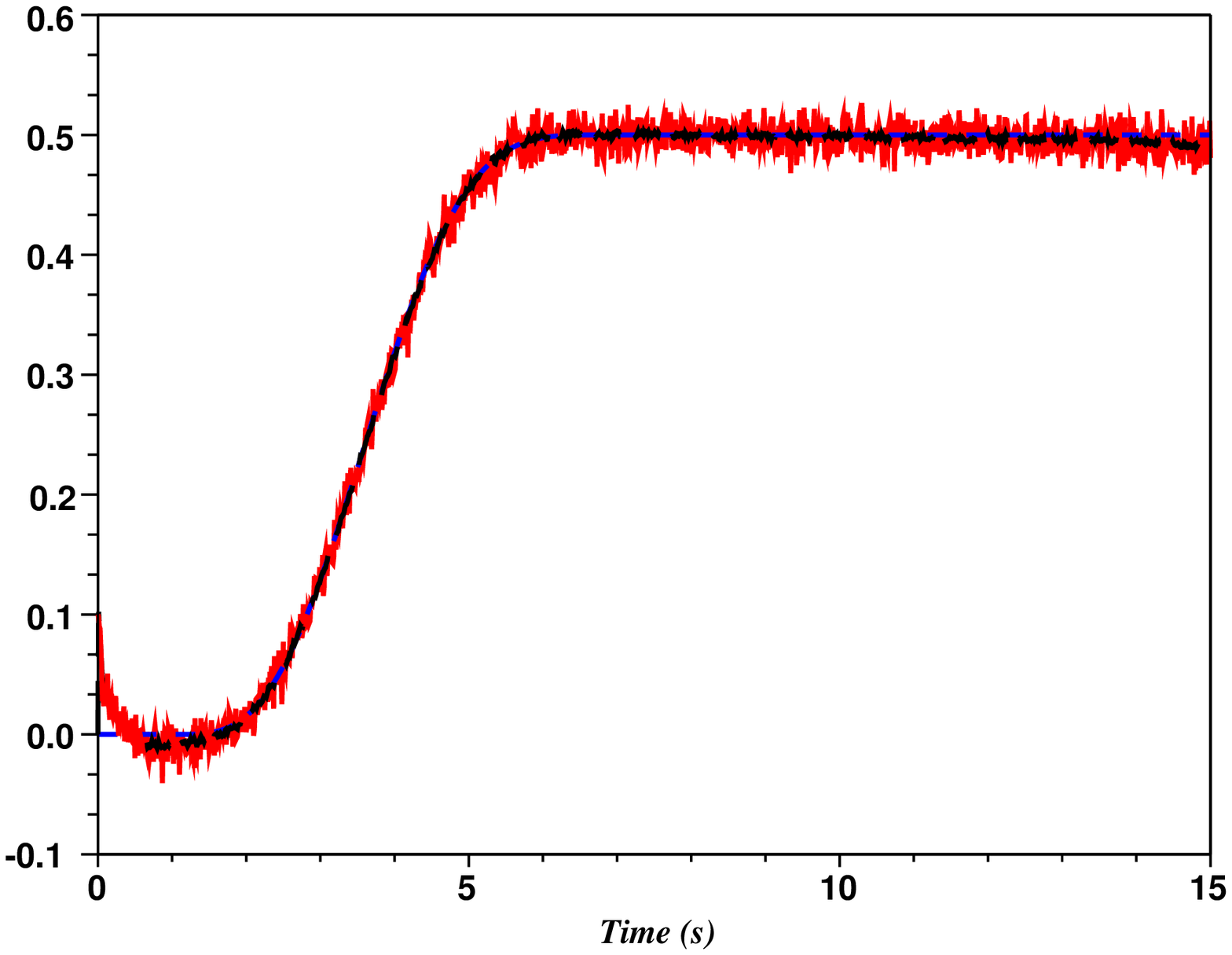}}
\caption{i-PI in case of power loss}\label{x1iPIPP}
\end{figure*}

Introduce now a fault accommodation via a control power loss $u_{\rm
Pert}=0.996^{t/h} \times u$, $t>4$, where the sampling time
$h=0.01s$. The i-PI behaves then much better (Figure \ref{x1iPIPP})
than the PI (Figure \ref{x1PIPP}). Note nevertheless a small
deviation of the i-PI controller when the power loss becomes quite
important (Figure \ref{x1iPIPP}-(b)).


%

\section{Conclusion}
\label{Conclusion} The above numerical simulations as well as many
existing experimentations (see \cite{esta,malo}, and
\cite{cifa-mines,choi,brest1,brest2,psa,edf,vil}) demonstrate that
intelligent PID controllers yield better performances than classic
ones. This is achieved moreover thanks to a quite straightforward
and natural gain tuning, which contrasts with the numerous complex
rules for classic PIDs. Those considerations as well as the results
of this communication imply therefore
\begin{itemize}
\item that classic PIDs might become obsolete,
\item a change of paradigm for
control engineering, and for its teaching (see, {\it e.g.},
\cite{esta,malo}, and \cite{riachy}).
\end{itemize}

\end{document}